# VARIANCE BOUNDING MARKOV CHAINS

By Gareth O. Roberts and Jeffrey S. Rosenthal[1]

*Lancaster University and University of Toronto*

We introduce a new property of Markov chains, called *variance bounding*. We prove that, for reversible chains at least, variance bounding is weaker than, but closely related to, geometric ergodicity. Furthermore, variance bounding is equivalent to the existence of usual central limit theorems for all $L^2$ functionals. Also, variance bounding (unlike geometric ergodicity) is preserved under the Peskun order. We close with some applications to Metropolis–Hastings algorithms.

**1. Introduction.** Markov chain Monte Carlo (MCMC) algorithms are widely used in statistics, physics, and computer science. Measures of how good an MCMC algorithm is include quantitative bounds on convergence to stationarity (e.g., [14, 15, 34, 35]), qualitative convergence rates such as geometric ergodicity (e.g., [21, 29, 32, 39, 40]), the existence of central limit theorems (e.g., [2, 3, 7, 10, 13, 21, 40]) and bounds on asymptotic variance of estimators (e.g., [7, 22, 41]).

In this paper we introduce a new notion, *variance bounding*. Roughly, a Markov chain is variance bounding if the asymptotic variances for functionals with unit stationary variance are uniformly bounded (precise definitions are given below). We shall show that, for reversible chains at least, variance bounding is implied by geometric ergodicity, and conversely, if $P$ is variance bounding, then $aI + (1-a)P$ is geometrically ergodic for all $0 < a < 1$. More importantly, we shall prove that a reversible Markov chain is variance bounding if and only if all $L^2$ functionals satisfy a usual central limit theorem, indicating that variance bounding is in some sense the "right" definition to use. We also prove that variance bounding is preserved under the Peskun partial ordering ([26, 40]) on Markov chains. Finally, applications to Metropolis–Hastings algorithms are presented.

Received October 2006; revised April 2007.
[1]Supported in part by NSERC of Canada.

**2. Variance bounding.** Given a Markov chain kernel $P$ on a state space $(\mathcal{X}, \mathcal{F})$ with unique stationary distribution $\pi(\cdot)$, we let $\{X_n\}$ follow the kernel $P$ in stationarity, so that $\mathbf{P}[X_n \in A] = \pi(A)$ for all $A \in \mathcal{F}$ and $n \in \mathbf{N} \cup \{0\}$, and also $\mathbf{P}[X_n \in A \mid X_0, \ldots, X_{n-1}] = P(X_{n-1}, A)$ for all $A \in \mathcal{F}$ and all $n \in \mathbf{N}$. For a functional $h : \mathcal{X} \to \mathbf{R}$ (assumed throughout to be measurable), the stationary variance is given by $\operatorname{Var}_\pi(h) = \mathbf{E}[(h(X_0) - \mathbf{E}[h(X_0)])^2]$, and the asymptotic variance is given by

$$
(1) \qquad \operatorname{Var}(h, P) = \lim_{n \to \infty} \frac{1}{n} \operatorname{Var}\left( \sum_{i=1}^{n} h(X_i) \right).
$$

If the Markov chain $P$ is to be used to estimate the stationary expected value of $h$ by $\frac{1}{n} \sum_{i=1}^{n} h(X_i)$, then $\operatorname{Var}(h, P)$ is a measure of the Monte Carlo uncertainty of the estimate. Thus, for MCMC algorithms, it is desirable to make $\operatorname{Var}(h, P)$ as small as possible (cf. [7, 22, 23, 40, 41]). This prompts the following definition.

DEFINITION. $P$ is *variance bounding* if there is $K < \infty$ such that $\operatorname{Var}(h, P) \leq K \operatorname{Var}_\pi(h)$ for all $h : \mathcal{X} \to \mathbf{R}$. Equivalently, $P$ is *variance bounding* if $\sup\{\operatorname{Var}(h, P); h : \mathcal{X} \to \mathbf{R}, \operatorname{Var}_\pi(h) = 1\} < \infty$.

Note that in the case where $\operatorname{Var}_\pi(h) = \infty$, the required inequality holds automatically for all $K$.

Variance bounding is a natural property, in that it offers some control over the asymptotic variances $\operatorname{Var}(h, P)$. We study its relation to more traditional MCMC properties below. For most of our results, we assume that $P$ is *reversible* with respect to $\pi(\cdot)$, that is, that

$$
(2) \qquad \int_{x \in A} \pi(dx) P(x, B) = \int_{x \in B} \pi(dx) P(x, A), \qquad A, B \in \mathcal{F}.
$$

It follows from [16] (see also [3]) that, for reversible chains and $L^2$ functionals, the limit in equation (1) always exists, though it may be infinite.

**3. Relation to geometric ergodicity.** Recall that a Markov chain kernel $P$ with stationary distribution $\pi(\cdot)$ is *geometrically ergodic* if there is $\rho < 1$ and $M : \mathcal{X} \to [0, \infty]$ $\pi$-a.e. finite [i.e., such that $\pi\{x \in \mathcal{X} : M(x) < \infty\} = 1$], such that $|P^n(x, A) - \pi(A)| \leq M(x)\rho^n$ for all $A \in \mathcal{F}$, $n \in \mathbf{N}$, and $x \in \mathcal{X}$. Geometric ergodicity is an often studied property (e.g., [21, 29, 32, 39, 40]), which leads to many useful results, such as central limit theorems (see next section).

However, geometric ergodicity is an overly strong notion in that it requires, among other things, that the Markov chain be *aperiodic*. Since estimates of functionals, and their variances $\operatorname{Var}(h, P)$, are essentially unaffected by periodicity considerations, it seems inappropriate to demand aperiodicity. And indeed, many Markov chains are variance bounding despite being periodic (e.g., the Markov chain $P_1$ in Example 9 below).

We now explore the relation between geometric ergodicity and variance bounding. We first show that, for reversible chains, variance bounding is strictly *weaker* than geometric ergodicity. (Proofs of all theorems are deferred until Section 7.)

THEOREM 1. *If $P$ is reversible and geometrically ergodic, then $P$ is variance bounding.*

Next, we show that $P$ is variance bounding if and only if any mixture of $P$ with the identity is geometrically ergodic. We write $I$ for the identity kernel, that is, the Markov chain which never moves, so that $I(x, \{x\}) = 1$ for all $x \in \mathcal{X}$.

THEOREM 2. *If $P$ is reversible, then the following are equivalent:*

(i) *$P$ is variance bounding.*
(ii) *$aI + (1-a)P$ is geometrically ergodic for all $0 < a < 1$.*
(iii) *$aI + (1-a)P$ is geometrically ergodic for some $0 \leq a < 1$.*

COROLLARY 3. *If $P$ is reversible, then for any fixed $0 \leq a < 1$, the following are equivalent:*

(i) *$P$ is variance bounding.*
(ii) *$aI + (1-a)P$ is variance bounding.*

Section 6 below contains some applications of Theorems 1 and 2. We next note that if $P$ has holding probabilities uniformly bounded away from 0, then variance bounding and geometrically ergodic are equivalent:

THEOREM 4. *If $P$ is reversible and $\inf_{x \in \mathcal{X}} P(x, \{x\}) > 0$, then $P$ is variance bounding if and only if $P$ is geometrically ergodic.*

As an application of Theorem 4, suppose $P$ represents a *random-walk Metropolis* or *systematic-scan Metropolis-within-Gibbs* algorithm on $\mathbf{R}^d$, with proposal increment densities positive in a neighborhood of 0, whose target density $t$ is $C^1$ with $\|\nabla \log t(x)\| \geq \delta > 0$ for all $x \in \mathcal{X}$. It then follows as in [33] that the rejection probabilities $P(x, \{x\})$ are uniformly bounded away from 0. Hence, by Theorem 4, variance bounding is equivalent to geometric ergodicity in this case.

Similarly, the two notions are equivalent if the operator $P$ is *positive*, that is, if $\mathbf{E}[f(X_0)f(X_1)] \geq 0$ for all measurable $f : \mathcal{X} \to \mathbf{R}$ when $\{X_n\}$ is in stationarity:

THEOREM 5. *If $P$ is reversible and positive, then $P$ is variance bounding if and only if $P$ is geometrically ergodic.*

As an application of Theorem 5, suppose $P$ represents a *data augmentation* algorithm, that is, the $x$-coordinate (only) of a two-variable Gibbs sampler. It follows from Lemmas 3-1 and 3-2 of [18] that $P$ is reversible and positive. Hence, by Theorem 5, variance bounding is equivalent to geometric ergodicity in this case as well. (See also [11].)

In particular, the *slice sampler* (e.g., [24, 25, 30]) can be viewed as the $x$-coordinate of a two-variable Gibbs sampler. (This holds for product slice samplers as well, since the multiple auxiliary variables are conditionally independent and can be regarded as a single auxiliary vector.) So, for any slice sampler, variance bounding is equivalent to geometric ergodicity. For example, it is known [30] that the slice sampler is geometrically ergodic whenever $Q'(y)y^{1+1/\alpha}$ is nonincreasing near 0, for some $\alpha > 1$, where $Q(y)$ is the measure of the set where the target density value is at least $y$. It follows immediately that the slice sampler is also variance bounding under these conditions.

In general, if $P$ is variance bounding, then a slight modification of $P$ is geometrically ergodic. Specifically, following [36], let $\overline{P}^n$ be the *binomial modification* of $P$, corresponding to doing an (independently chosen) random number $B_n$ of steps from $P$, where $B_n \sim \text{Binomial}(2n, 1/2)$. Thus, $\overline{P}^n = 2^{-2n} \sum_{i=0}^{2n} \binom{2n}{i} P^i$. Call $\overline{P}$ geometrically ergodic if, as usual, there is $\rho < 1$ and $\pi$-a.e. finite $M : \mathcal{X} \to [0, \infty]$ such that $|\overline{P}^n(x, A) - \pi(A)| \leq M(x)\rho^n$ for all $A \in \mathcal{F}$, $n \in \mathbf{N}$, and $x \in \mathcal{X}$. Then we have the following.

THEOREM 6. *If $P$ is reversible, then $P$ is variance bounding if and only if $\overline{P}$ is geometrically ergodic.*

REMARK. The stationary processes literature (e.g., [2, 12, 13]) defines many other mixing conditions, such as $\alpha$-mixing, $\beta$-mixing, $\rho$-mixing, $\phi$-mixing, etc. These conditions are related to usual Markov chain ergodicity conditions, for example, $\phi$-mixing is equivalent to uniform ergodicity, exponentially-fast $\beta$-mixing is equivalent to geometric ergodicity, $\alpha$-mixing is implied by Harris ergodicity, etc. However, none of these mixing conditions is implied by variance bounding, since the mixing conditions all require ergodicity, whereas periodic (and therefore nonergodic) chains can still be variance bounding.

**4. Relation to central limit theorems.** An important issue in MCMC is the existence of central limit theorems (e.g., [2, 3, 7, 10, 13, 40]). Where central limit theorems are known to hold, they underpin practical MCMC strategies for Monte Carlo error assessment (see, e.g., [8]).

Say that a functional $h:\mathcal{X} \to \mathbf{R}$ with $\pi(|h|) < \infty$ [where $\pi(f) = \int_{\mathcal{X}} f(x) \times \pi(dx)$] satisfies a usual central limit theorem (CLT) for a Markov chain $P$ if, as $n \to \infty$, the distribution of $n^{-1/2} \sum_{i=1}^{n} [h(X_i) - \pi(h)]$ converges weakly to $N(0, v)$, where (with $\{X_n\}$ in stationarity)

$$(3) \qquad v = \mathrm{Var}(X_0) + 2 \sum_{i=1}^{\infty} \mathrm{Cov}(X_0, X_i) < \infty.$$

(We say "usual" to distinguish this convergence from, e.g., convergence to other distributions, or other normalizations besides $n^{-1/2}$; see also [3, 7, 27, 37].)

It is known ([12], Theorem 18.5.3; see also [3], [10]) that if $P$ is geometrically ergodic, then $h$ satisfies a usual CLT, provided $\pi(|h|^{2+\delta}) < \infty$ for some $\delta > 0$. It was proven in [28], following [16], that if $P$ is geometrically ergodic and reversible, then $h$ satisfies a usual CLT whenever $\pi(h^2) < \infty$. However, geometric ergodicity is an overly strong assumption; for example, periodic Markov chains can never be geometrically ergodic but they can still satisfy CLTs.

The following theorem shows that, for reversible Markov chains, variance bounding is the "right" definition for CLTs, that is, variance bounding (unlike geometric ergodicity) is the weakest property which still guarantees usual CLTs for all $L^2$ functionals. (We assume the stationary distribution for $P$ is unique, to avoid degenerate cases where the state space breaks up into multiple closed subsets.)

THEOREM 7. *If $P$ is reversible, with unique stationary distribution $\pi(\cdot)$, then $P$ is variance bounding if and only if every $h:\mathcal{X} \to \mathbf{R}$ with $\pi(h^2) < \infty$ satisfies a usual CLT for $P$.*

REMARK. There are other results available (see, e.g., [3] and the references therein) which guarantee CLTs for *specific* functionals, rather than for all $L^2$ functionals. However, often MCMC is used to generate samples from $\pi(\cdot)$ before it is decided which functionals are of statistical interest. Thus, we find that it is most useful having results like Theorem 7 which apply to all $L^2$ functionals simultaneously.

**5. Relation to the Peskun ordering.** The following partial order on Markov chain kernels was introduced by Peskun [26] for finite state spaces, and later by Tierney [40] for general state spaces.

DEFINITION. Let $P_1$ and $P_2$ be two Markov chain kernels on $(\mathcal{X}, \mathcal{F})$, both having invariant probability measure $\pi$. Then $P_1$ *dominates* $P_2$ *off the diagonal*, written $P_1 \succeq P_2$, if $P_1(x, A) \geq P_2(x, A)$ for all $x \in \mathcal{X}$ and $A \in \mathcal{F}$ with $x \notin A$.

It was proved by Peskun [26] for finite state spaces, and then by Tierney [41] (see also [22, 23]) for general state spaces, that if $P_1 \succeq P_2$, and $P_1$ and $P_2$ are reversible with respect to the same $\pi(\cdot)$, then $\text{Var}(h, P_1) \leq \text{Var}(h, P_2)$ for all $h : \mathcal{X} \to \mathbf{R}$. That is, $P_1$ is "better" than $P_2$, in the sense of being uniformly more efficient for estimating expectations of functionals. Thus, it seems reasonable that any Markov chain property designed to indicate good estimation should be preserved under the Peskun ordering. For the variance bounding property, that is indeed the case:

THEOREM 8. *If $P_1$ and $P_2$ are both reversible with respect to $\pi(\cdot)$, and $P_1 \succeq P_2$, and $P_2$ is variance bounding, then $P_1$ is variance bounding.*

On the other hand, the corresponding property for geometric ergodicity does not hold, indicating another advantage of variance bounding over geometric ergodicity:

EXAMPLE 9. Let $\mathcal{X} = \mathbf{Z}$ with $\pi(m) = 2^{-|m|}/3$. Define $P_1$ by $P_1(x, x-1) = 2/3$ and $P_1(x, x+1) = 1/3$ for $x > 0$, and $P_1(x, x-1) = 1/3$ and $P_1(x, x+1) = 2/3$ for $x < 0$, and $P_1(0, -1) = P_1(0, 1) = 1/2$. Also, let $P_2$ be the Metropolis algorithm for $\pi(\cdot)$ with proposal distribution $Q(x, x+1) = Q(x, x-1) = 1/2$. [Thus, $P_2(x, x+1) = P_2(x, x) = 1/4$ and $P_2(x, x-1) = 1/2$ for $x > 0$; $P_2(x, x-1) = P_2(x, x) = 1/4$ and $P_2(x, x+1) = 1/2$ for $x < 0$; and $P_2(0, -1) = P_2(0, 1) = 1/4$ and $P_2(0, 0) = 1/2$.] Then both $P_1$ and $P_2$ are reversible with respect to $\pi(\cdot)$, and also $P_1 \succeq P_2$. Furthermore, it follows as in Mengersen and Tweedie [19] that $P_2$ is geometrically ergodic, and hence variance bounding by Theorem 1. On the other hand, $P_1$ is periodic, and hence cannot be geometrically ergodic, even though $P_1 \succeq P_2$. (Of course, $P_1$ is still variance bounding, by Theorem 8.)

**6. Application to Metropolis–Hastings algorithms.** We now consider Metropolis–Hastings algorithms ([9, 20]). We define a slight generalization, as follows. Given a reference measure $\nu(\cdot)$ on $\mathcal{X}$, with respect to which $\pi(dx) = t(x)\nu(dx)$, and a nonnegative (measurable) function $q : \mathcal{X} \times \mathcal{X} \to \mathbf{R}$ with $\int_{\mathcal{X}} q(x, y)\nu(dy) \leq 1$ for all $x \in \mathcal{X}$, the *sub-Metropolis–Hastings algorithm* is the algorithm with transition kernel

$$M_q(x, dy) = \alpha(x, y)q(x, y)\nu(dy) + r(x)\delta_x(dy),$$

where $\alpha(x, y) = \min(1, \frac{t(y)q(y,x)}{t(x)q(x,y)})$, and $r(x) = 1 - \int_{\mathcal{X}} \alpha(x, y)q(x, y)\nu(dy) \geq 0$.

By construction, this algorithm is reversible with respect to $\pi(\cdot)$. It may be described as follows. With probability $\int_{\mathcal{X}} q(x, y)\nu(dy)$, it performs the usual Metropolis–Hastings algorithm with proposal density $q(x, y)/\int_{\mathcal{X}} q(x, y)\nu(dy)$. Otherwise, with probability $1 - \int_{\mathcal{X}} q(x, y)\nu(dy)$, it stays at its current state. If $\int_{\mathcal{X}} q(x, y)\nu(dy) = 1$, then $M_q$ is the usual Metropolis–Hastings algorithm.

By direct inspection, noting that $\alpha(x,y)q(x,y) = \min(q(x,y), \frac{t(y)}{t(x)}q(y,x))$, we see the following:

PROPOSITION 10. *For fixed $\nu(\cdot)$ and $t$, if $q_1(x,y) \geq q_2(x,y)$ for all $x,y \in \mathcal{X}$ with $x \neq y$, then $M_{q_1} \succeq M_{q_2}$. (Hence, by Theorem 8, if $M_{q_2}$ is variance bounding, then so is $M_{q_1}$.)*

Now, suppose $M_{q_2}$ is variance bounding, and that $q_1(x,y) \geq cq_2(x,y)$ for all $x \neq y$, for some $c > 0$. We can assume [by replacing $c$ with $\min(c,1)$ if necessary] that $c \leq 1$. Then $M_{cq_2} = cM_{q_2} + (1-c)I$. Hence, by Corollary 3 (with $a = 1-c$), $M_{cq_2}$ is also variance bounding. It then follows from Proposition 10 that $M_{q_1}$ is also variance bounding. We conclude:

COROLLARY 11. *If $q_1(x,y) \geq cq_2(x,y)$ for all $x,y \in \mathcal{X}$ with $x \neq y$, for some $c > 0$, and if $M_{q_2}$ is variance bounding, then $M_{q_1}$ is variance bounding.*

Example 9 above shows that the analogous statement to Corollary 11 for geometric ergodicity does not hold.

To continue, call a (measurable) function $s: \mathcal{X} \to [0, \infty)$ *MT-good* if it is symmetric, positive and continuous, with exponentially bounded tails, and with $\int_{-\infty}^{\infty} s(u)\,du = 1$. Then a result of Mengersen and Tweedie [19] (see also [32] for higher-dimensional analogs) says that a random-walk Metropolis algorithm on $\mathcal{X} = \mathbf{R}$, with proposal density $q(x,y) = s(y-x)$ for some MT-good $s$, is geometrically ergodic provided the target density has exponentially bounded tails. This is a very impressive result, but with the severe restriction that the proposal increments must correspond to a symmetric random walk. To improve this, we make the following definition.

DEFINITION. A proposal density function $q: \mathcal{X} \times \mathcal{X} \to \mathbf{R}$ is a *uniformly minorized increment distribution* (*UMID*) if there is $c > 0$ and MT-good $s: \mathcal{X} \to [0, \infty)$ such that $q(x,y) \geq cs(y-x)$ for all $x, y \in \mathcal{X}$.

Combining Theorem 1 and Corollary 11 with the result of [19] immediately gives the following:

COROLLARY 12. *Let $t$ be a target density with exponentially bounded tails, and let $q$ be a UMID proposal density function. Then $M_q$ is variance bounding.*

Note that in Corollary 12, we do not need to assume that $s$ has exponentially bounded tails, since if not then we can simply replace $s(x)$ by $\min(s(x), e^{-|x|})$ without affecting the conclusion. Note also that Corollary 12 does not require the proposal density $q$ to be symmetric, nor to correspond to a random walk. (Similar generalizations are also available for the multi-dimensional case, as in [32].)

As one application of Corollary 12, consider a Langevin (MALA) algorithm (see [33]), with proposal density given by $Q(x, \cdot) = N(x + \frac{1}{2}\delta \nabla \log t(x), \delta^2)$ for some $\delta > 0$. Now, if the target density $t$ is $C^1$ with tails that are precisely exponential, then $\nabla \log t(x)$ is a bounded function of $x \in \mathcal{X}$, and it follows easily that $q$ is UMID. We conclude:

COROLLARY 13. *A Langevin algorithm for a $C^1$ target density on $\mathcal{X} = \mathbf{R}$ with exponential tails is variance bounding.*

As a final application, we consider a Metropolis–Hastings algorithm for a density $t$ supported on $(0, \infty)$, with proposal distribution given by $Q(x, \cdot) = N(x, x^b)$ for some fixed $b > 0$. That is, the variance of the proposal increment depends on the current state $x \in \mathcal{X}$. (Related models were considered in [31].)

If $b > 2$, then as $x \to \infty$, the proposal values will be farther and farther out in the tails, so $\lim_{x \to \infty} P(x, \{x\}) = 1$. It follows as in [32], or by a simple capacitance argument (e.g., [17]), that the resulting Markov chain is neither geometrically ergodic nor variance bounding. So, we do not consider that case further here. (On the other hand, numerical simulations related to [31] indicate that if $t$ is, e.g., a Cauchy distribution, then values $b \approx 2.7$ may give fastest numerical convergence, which is a separate but related issue.)

If $b = 2$, then the distribution $Q(x, \cdot)$ equals the distribution of $x + xZ$, where $Z \sim N(0, 1)$. Taking logarithms (cf. [31]) gives rise to an equivalent chain which is an ordinary random-walk Metropolis algorithm, with modified target density $\tilde{t}(y) = e^y t(e^y)$, and with increment density $f(u)$ equal to the density of $\log(1 + Z)$ where $Z \sim N(0, 1)$. This increment distribution is clearly UMID; indeed, we can simply let $cs(u) = \min(f(u), f(-u))$. Hence, by Corollary 12, the transformed chain—and hence, also the original chain—is variance bounding, provided that $\tilde{t}$ has exponentially bounded tails.

Finally, suppose that $0 < b < 2$. Then $Q(x, \cdot)$ is the distribution of $x + x^{b/2}Z$, where $Z \sim N(0, 1)$. Instead of logarithms, consider the transformation $X \mapsto X^a$, where $a = 1 - b/2$ (so $0 < a < 1$). Then the proposal increment from $x \in \mathcal{X}$ transforms from $x^{b/2}Z$ to $W = h(Z) \equiv [x + x^{b/2}Z]^a - x^a = x^a(1 + x^{-a}Z)^a - x^a$. Inverting this, $Z = h^{-1}(W) = x^a((1 + Wx^{-a})^{1/a} - 1)$. Now, the density of $Z$ is $g(z) = (2\pi)^{-1/2} e^{-z^2/2}$. Hence, for the transformed chain, the proposal increment $W$ has density

$$\frac{g(h^{-1}(w))}{(dw/dz)} = \frac{e^{-h^{-1}(w)^2/2}}{\sqrt{2\pi} a(1 + x^{-a}z)^{a-1}}.$$

We compute that, as $x \to \infty$, this expression converges to $(2\pi)^{-1/2}a^{-1}e^{-(w/a)^2/2}$, that is, to the density function of the $N(0, a^2)$ distribution. [Intuitively, this is because $(\frac{d}{dx}x^a)^2 x^b = a^2 x^{2a-2} x^b = a^2$ is constant, so the increment variance of the transformed chain is approximately stabilized.] Hence, for large enough $x$, and thus for all $x$ by positivity and continuity, the proposal density is UMID. Therefore, by Corollary 12, the transformed and original chains are variance bounding, provided that the transformed target density has exponentially bounded tails.

**7. Spectra and theorem proofs.** We now proceed to the proofs of the theorems. We begin by recalling some standard notation. Let $P$ be a Markov chain kernel with stationary distribution $\pi(\cdot)$ on a state space $(\mathcal{X}, \mathcal{F})$. For measurable $f, g : \mathcal{X} \to \mathbf{R}$, write $\langle f, g \rangle = \int_{\mathcal{X}} f(x)g(x)\pi(dx)$, and $\|f\| = \langle f, f \rangle^{1/2}$. Let $L_0^2(\pi) = \{f : \mathcal{X} \to \mathbf{R}; \pi(f) = 0, \pi(f^2) < \infty\}$, and regard $P$ as an operator acting on $L_0^2(\pi)$, by $(Pf)(x) = \int_{\mathcal{X}} f(y)P(x, dy)$. Write $\sigma(P)$ for the spectrum of the operator $P$ acting on $L_0^2(\pi)$ (see, e.g., [4, 36]). If $P$ is a reversible Markov chain, then $P$ is a self-adjoint operator with respect to $\langle \cdot, \cdot \rangle$, and also $\sigma(P) \subseteq [-1, 1]$ (cf. [1, 7]). Theorem 2 of [28] says that if $P$ is reversible, then $P$ is geometrically ergodic if and only if there is $r < 1$ with $\sigma(P) \subseteq [-r, r]$.

We have the following.

THEOREM 14. *If $P$ is reversible, then $P$ is variance bounding if and only if $\sup(\sigma(P)) < 1$.*

PROOF. Suppose first that $\sup(\sigma(P)) \equiv \Lambda < 1$. Then by Proposition 1 of [36], $\text{Var}(h, P) \leq 2(1-\Lambda)^{-1} \text{Var}_\pi(h)$ for all $h : \mathcal{X} \to \mathbf{R}$. Hence, $P$ is variance bounding with constant $K = 2(1-\Lambda)^{-1} < \infty$.

Conversely, suppose $\sup(\sigma(P)) = 1$. Let $\mathcal{E}$ be the spectral measure for $P$ (see, e.g., [4, 7, 28, 36, 38]), and let $r < 1$. Then $\mathcal{E}((r, 1])$ is nonzero, so there is $h \in L_0^2(\pi)$ in range of $\mathcal{E}((r, 1])$. It follows similarly to Proposition 1 of [36] (cf. [7, 16]) that $\text{Var}(h, P) \geq \frac{1}{1-r} \text{Var}_\pi(h)$. Since this holds for any $r < 1$, it follows that $\sup_{h \in \mathcal{L}_0^2(\pi)} [\text{Var}(h, P)/\text{Var}_\pi(h)] \geq \sup_{r<1} \frac{1}{1-r} = \infty$. Hence, $P$ is not variance bounding. □

PROOF OF THEOREM 1. If $P$ is reversible and geometrically ergodic, then there is $r < 1$ with $\sigma(P) \subseteq [-r, r]$. In particular, $\sup(\sigma(P)) \leq r < 1$, so $P$ is variance bounding by Theorem 14. □

PROOF OF COROLLARY 2. (i) $\Longrightarrow$ (ii): Suppose $P$ is variance bounding, and $0 < a < 1$. Then by Theorem 14, $\sup(\sigma(P)) < 1$, that is, there is $c < 1$ with $\sigma(P) \subseteq [-1, c]$. On the other hand,

$$\sigma(aI + (1-a)P) = \{\lambda \in \mathbf{R} \text{ s.t. } (aI + (1-a)P - \lambda I) \text{ is not invertible}\}$$

$$= \left\{\lambda \in \mathbf{R} \text{ s.t. } (1-a)\left(P - \frac{\lambda - a}{1-a}I\right) \text{ is not invertible}\right\}$$

$$= \left\{\lambda \in \mathbf{R} \text{ s.t. } \frac{\lambda - a}{1-a} \in \sigma(P)\right\}$$

$$= \{(a + (1-a)y) \text{ s.t. } y \in \sigma(P)\},$$

where the last equality follows by solving for $\lambda$ in the equation $y = \frac{\lambda - a}{1-a}$. Hence, since $\sigma(P) \subseteq [-1, c]$, it follows that

$$\sigma(aI + (1-a)P) \subseteq [a + (1-a)(-1), a + (1-a)c]$$
$$= [2a - 1, a + (1-a)c] \subseteq [-r, r],$$

where $r = \max(|2a - 1|, a + (1-a)c) < 1$. Hence, by Theorem 2 of [28], $aI + (1-a)P$ is geometrically ergodic.

(ii) $\implies$ (iii): Immediate.

(iii) $\implies$ (i): If $aI + (1-a)P$ is geometrically ergodic, then there is $r < 1$ with $\sigma(aI + (1-a)P) \subseteq [-r, r]$. But from the above,

$$\sigma(aI + (1-a)P) = \left\{\lambda \in \mathbf{R} \text{ s.t. } \frac{\lambda - a}{1-a} \in \sigma(P)\right\},$$

so it follows that $\sigma(P) \subseteq [\frac{-r-a}{1-a}, \frac{r-a}{1-a}]$. In particular, $\sup(\sigma(P)) \leq \frac{r-a}{1-a} < 1$, so $P$ is variance bounding. $\square$

PROOF OF THEOREM 3. We see from the proof of Theorem 2 that

$$\sup(\sigma(aI + (1-a)P)) = a + (1-a)\sup(\sigma(P)).$$

It follows that for $0 \leq a < 1$, $\sup(\sigma(aI + (1-a)P)) < 1$ if and only if $\sup(\sigma(P)) < 1$. The result then follows from Theorem 14. $\square$

PROOF OF THEOREM 4. If $P$ is reversible and geometrically ergodic, then $P$ is variance bounding by Theorem 1. Conversely, suppose $P$ is reversible and variance bounding, with $\delta \equiv \inf_{x \in \mathcal{X}} P(x, \{x\}) > 0$. Let $S(x, A) = (1 - \delta)^{-1}(P(x, A) - \delta \mathbf{1}_{x \in A})$. Then $S$ is another Markov chain kernel on $\mathcal{X}$, and $P = \delta I + (1 - \delta)S$. It follows that $\inf \sigma(P) \geq \delta + (1 - \delta)(-1) = 2\delta - 1 > -1$. Also $\sup \sigma(P) < 1$ by Theorem 14. Hence, there is $r < 1$ with $\sigma(P) \subseteq [-r, r]$, so $P$ is geometrically ergodic. $\square$

PROOF OF THEOREM 5. Note that $\mathbf{E}[f(X_0)f(X_1)] = \langle f, Pf \rangle$, so positivity is equivalent to $\langle f, Pf \rangle \geq 0$ for all $f \in L_0^2(\pi)$. This implies that $\lambda \geq 0$ for all $\lambda \in \sigma(P)$. Hence, using Theorem 14,

$$P \text{ is geometrically ergodic} \iff \sup\{|\lambda| : \lambda \in \sigma(P)\} < 1$$
$$\iff \sup\{\lambda : \lambda \in \sigma(P)\} < 1$$
$$\iff P \text{ is variance bounding.} \qquad \square$$

PROOF OF THEOREM 6. Note that we can write $\overline{P}^n = [\frac{1}{2}(I + P)]^n$. Hence, the result follows immediately from Theorem 2 (with $a = 1/2$). $\square$

PROOF OF THEOREM 7. If $P$ is variance bounding, then $\Lambda \equiv \sup(\sigma(P)) < 1$. Let $\mathcal{E}$ be the spectral measure for $P$, and let $\mathcal{E}_h$ be the induced measure defined by

$$\mathcal{E}_h(S) = \int_X h(x)(\mathcal{E}(S)h)(x)\pi(dx).$$

Then it follows (cf. [7]) that

$$\sigma^2 \equiv \int_{-1}^{1} \frac{1+\lambda}{1-\lambda} \mathcal{E}_h(d\lambda) \leq \frac{1+\Lambda}{1-\Lambda} < \infty.$$

It then follows from Kipnis and Varadhan [16] (see also [3]) that $h$ satisfies a usual CLT for $P$.

Conversely, if $P$ is not variance bounding, then $\Lambda = 1$. It follows as in the proof of Theorem 14 that $\mathcal{E}((r,1])$ is nonzero for every $r < 1$. Since $P$ has unique stationary distribution, $1 \notin \sigma(P)$, so there must be infinitely many $m \in \mathbf{N}$ such that $\mathcal{E}((1 - 2^{-m}, 1 - 2^{-m-1}])$ is nonzero. Let $m_1 < m_2 < \cdots$ (so $m_i \geq i$) with $\mathcal{E}((1 - 2^{-m_i}, 1 - 2^{-m_i-1}])$ bee nonzero. Let $g_i \in L_0^2(\pi)$ be in the range of $\mathcal{E}((1 - 2^{-m_i}, 1 - 2^{-m_i-1}])$, with $\|g_i\| = 1$. Then spectral theory implies that the $\{g_i\}$ are orthonormal, and furthermore, $\text{Cov}(g_i, Pg_i) = \langle g_i, Pg_i \rangle \geq 1 - 2^{-m_i}$. Finally, let $h = \sum_{i=1}^{\infty} 2^{-i/2} g_i$. Then by orthonormality,

$$\text{Var}_\pi(h) = \|h\|^2 = \sum_{i=1}^{\infty} (2^{-i/2})^2 = 1 < \infty.$$

On the other hand, with $\{X_n\}$ in stationarity, again using orthonormality,

$$\text{Cov}(h(X_0), h(X_n)) = \sum_{i=1}^{\infty} 2^{-i} \text{Cov}(g_i, P^n g_i)$$
$$\geq \sum_{i=1}^{\infty} 2^{-i}(1 - 2^{-m_i})^n$$
$$\geq \sum_{i=1}^{\infty} 2^{-i}(1 - 2^{-i})^n.$$

Hence,

$$\sum_{n=0}^{\infty} \text{Cov}(h(X_0), h(X_n)) \geq \sum_{i=1}^{\infty} 2^{-i} \sum_{n=0}^{\infty} (1 - 2^{-i})^n$$

$$= \sum_{i=1}^{\infty} 2^{-i} [1 - (1 - 2^{-i})]^{-1}$$

$$= \sum_{i=1}^{\infty} 2^{-i} 2^i = \sum_{i=1}^{\infty} (1) = \infty.$$

It follows that $v$ in (3) is infinite, so $h$ does not satisfy a usual CLT for $P$. □

PROOF OF THEOREM 8. Lemma 3 of Tierney [41] says that since $P_1 \succeq P_2$, therefore $P_2 - P_1$ is a positive operator. It follows that $\sup(\sigma(P_2)) \geq \sup(\sigma(P_1))$. Hence, using Theorem 14 twice, if $P_2$ is variance bounding, then $\sup(\sigma(P_2)) < 1$, so $\sup(\sigma(P_1)) < 1$, so $P_1$ is variance bounding. [Alternatively, by Theorem 4 of [41], $\text{Var}(h, P_1) \leq \text{Var}(h, P_2) \leq K\text{Var}_\pi(h)$.] □

REMARK 1. The above theorems have all been proven for reversible chains only. However, it seems likely that analogs of some of them (e.g. Theorem 1) carry over in some form to nonreversible chains, about which various facts about convergence are known (see, e.g., [5, 6, 18, 23]). We leave this as an open problem for future work.

**8. Summary.** This paper defined a Markov chain to be variance bounding if the asymptotic variances for functionals with unit stationary variance are uniformly bounded. For reversible chains, we proved that this property is weaker than geometric ergodicity, but equivalent to $aI + (1-a)P$ being geometrically ergodic for all $0 < a < 1$. Furthermore, in contrast to geometric ergodicity, the variance bounding property: allows for periodicity; is equivalent to all $L^2$ functionals satisfying a usual central limit theorem; and is preserved under the Peskun [26] partial ordering on Markov chains. We also presented some applications to Metropolis–Hastings MCMC algorithms, and showed how variance bounding could apply more easily and more generally than geometric ergodicity.

Overall, we view these results as indicating that as a property to use in the study of MCMC algorithms, variance bounding is similar to, but more convenient than, geometric ergodicity. We hope that the notion of variance bounding can be used to further understand Markov chains and MCMC algorithms in other contexts.

**Acknowledgments.** We thank the anonymous referee for many very helpful comments.

DEPARTMENT OF MATHEMATICS AND STATISTICS
FYLDE COLLEGE
LANCASTER UNIVERSITY
LANCASTER, LA1 4YF
UNITED KINGDOM
E-MAIL: [g.o.roberts@lancaster.ac.uk](g.o.roberts@lancaster.ac.uk)

DEPARTMENT OF STATISTICS
UNIVERSITY OF TORONTO
TORONTO, ONTARIO
CANADA M5S 3G3
E-MAIL: [jeff@math.toronto.edu](jeff@math.toronto.edu)
URL: [http://probability.ca/jeff/](http://probability.ca/jeff/)